\documentclass[a4paper]{jpconf}
\usepackage{graphicx}
\usepackage[latin1]{inputenc}
\usepackage{amssymb}
\usepackage{amsmath}
\usepackage{latexsym}
\usepackage{cite}
\usepackage[frame,ps,matrix,arrow,curve,rotate,all,2cell,tips]{xy}
\usepackage{amssymb}
\usepackage{amsfonts}
\usepackage[normalem]{ulem}
\usepackage{amscd}
\usepackage{xcolor}
\usepackage{hyperref}
\usepackage{amstext,amsmath,amssymb,amsfonts}
\usepackage[latin1]{inputenc}
\usepackage{epsfig}
\usepackage{hyperref}
\newtheorem{theorem}{Theorem}

\newtheorem{corollary}[theorem]{Corollary}

\newtheorem{definition}[theorem]{Definition}
\newtheorem{example}[theorem]{Example}

\newtheorem{proposition}[theorem]{Proposition}
\newtheorem{remark}[theorem]{Remark}

\newcommand{\A}{\mathcal{A}}

\newcommand{\B}{\mathcal{B}}

\newcommand{\beq}{\begin{eqnarray}}
\newcommand{\eeq}{\end{eqnarray}}
\newcommand{\beqs}{\begin{eqnarray*}}
	\newcommand{\eeqs}{\end{eqnarray*}}
\newcommand{\bpro}{\begin{pro}}
	\newcommand{\epro}{\end{pro}}
\newcommand{\blem}{\begin{lem}}
	\newcommand{\elem}{\end{lem}}
\newcommand{\bdfn}{\begin{dfn}}
	\newcommand{\edfn}{\end{dfn}}
\newcommand{\bcor}{\begin{cor}}
	\newcommand{\ecor}{\end{cor}}
\newcommand{\bthm}{\begin{thm}}
	\newcommand{\ethm}{\end{thm}}
\newcommand{\bex}{\begin{ex}}
	\newcommand{\eex}{\end{ex}}
\newcommand{\brmk}{\begin{rmk}}
	\newcommand{\ermk}{\end{rmk}}
\newcommand{\bpr}{\begin{pr}}
	\newcommand{\epr}{\end{pr}}
\newcommand{\benum}{\begin{enumerate}} 
	\newcommand{\eenum}{\end{enumerate}}
\newcommand{\bitem}{\begin{itemize}}
	\newcommand{\eitem}{\end{itemize}}

\newcommand{\cqfd}{\hfill{\square}}
\numberwithin{equation}{subsection}
\numberwithin{table}{section}
\numberwithin{theorem}{section}
\DeclareMathOperator{\id}{id}

\begin{document}
\title{Zinbiel algebras and bialgebras: main properties and related  algebraic structures}

\author{Mahouton Norbert Hounkonnou}

\address{University of Abomey-Calavi,
	International Chair in Mathematical Physics and Applications,
	ICMPA-UNESCO Chair, 072 BP 50, Cotonou, Rep. of Benin}

\ead{norbert.hounkonnou@cipma.uac.bj, with copy to hounkonnou@yahoo.fr}

\author{Mafoya Landry Dassoundo}
\address{University of Abomey-Calavi,
	International Chair in Mathematical Physics and Applications,
	ICMPA-UNESCO Chair, 072 BP 50, Cotonou, Rep. of Benin}
\ead{mafoya.dassoundo@cipma.uac.bj}
\begin{abstract}
This work provides a characterization of  left  and right Zinbiel algebras.
Basic identities are established and discussed, showing that Zinbiel algebras are center-symmetric,  
and therefore Lie-admissible algebras. Their
bimodules  are given, and used to build a Zinbiel algebra structure
on the direct sum of the underlying vector space and a finite-dimensional vector space. In addition, their matched pair  is built,
and related to the matched pair of their sub-adjacent Lie algebras.
Besides, Zinbiel coalgebras are introduced, and linked to their underlying Lie coalgebras and
coassociative coalgebras. 
Moreover, the related Manin triple is defined, and used to characterize Zinbiel bialgebras, and their equivalence to the 
associated matched pair.
\end{abstract}

\section{Introduction}
	One of the important  subclasses  
	of nonassociative algebras
is that of Lie algebras. 
Their properties 
lead to
some generalizations of associative algebras such as 
Lie-admissible algebras\cite{Santilli},  left-symmetric algebras\cite{Bai_LSA}, flexible algebras,
alternative algebras, Malcev algebras,  Leibniz algebras, etc.

Leibniz algebra is an algebra $(\mathcal{A}, [,])$ such that the following identity
	is satisfied:
	\begin{eqnarray}
		[x,[y,z]]=[[x,y],z]-[[x,z], y], \; \forall x,y,z\in \mathcal{A}.
	\end{eqnarray}
 The category of Zinbiel algebras is  Koszul dual to the category of Leibniz algebras
	in the sense of J-L. Loday\cite{Loday}, 		
	where Zinbiel algebras were called 
	Leibniz algebras. 
	They appear in various domains in mathematics and physics, e.g. ,   in
		the theory of nonlinear geometric  control\cite{Kawski} characterized by the following system:
\begin{eqnarray}
	\left\lbrace
	\begin{array}{ccc}
	\dot{x}_1&=& u_1, \qquad |u_1|\leqslant 1\\
	\dot{x}_2&=& u_2, \qquad |u_2|\leqslant 1\\
	\dot{x}_3&=& x_1u_2-u_2 x_1. 
	\end{array}
	\right.
\end{eqnarray}
This system is completely controllable, i.e.  every point can be reached from every other point, since the
	distribution spanned by the vector fields $f_1=\partial_{x_1}-x_2 \partial_{x_3}$ and $f_2=\partial_{x_1}-x_3\partial_{x_3}$
	is nonintegrable due
	to the nonvanishing commutator $[f_1, f_2]=2\partial_{x_3}$.
	Performing the global
	coordinate change
	$y_1 = x_1,\; y_2 = x_2,$  $
		y_3 = \frac{1}{2}(x_3 + x_1 x_2 )$ 
	transforms  the initial system into a less symmetric but
	simpler form:
\begin{eqnarray}
	\left\lbrace
	\begin{array}{ccc}
	\dot{y}_1=&u_1,& \qquad |u_1|\leqslant 1\\
	\dot{y}_2=&u_2,& \qquad |u_2|\leqslant 1\\
	\dot{y}_3=&y_1u_2.&
	\end{array}
	\right.
\end{eqnarray}
Now setting
\begin{eqnarray}
	u_3(t)=(u_1\star u_2)(t):=\left(\int_{0}^{t}u_1(s)ds\right)u_2(t)
\end{eqnarray}
yields the Zinbiel identity 
\begin{eqnarray}
	(u_1\star(u_2\star u_3 )) (t)=((u_1\star u_2 + u_2\star u_1)\star u_3)(t).
\end{eqnarray}

Besides, Zinbiel algebras under q-commutator  given by  $x\circ_q y=x\circ y+qy\circ x$, where $q\in \mathbb{C}$, were investigated in \cite{Dzhumadil'daev}. 

In this work, we provide a characterization of left and right Zinbiel algebras. We establish  and discussed their
basic identities, showing that Zinbiel algebras are center-symmetric,
and, therefore,  Lie-admissible algebras. Furthermore, the bimodules are defined and used to
build a Zinbiel algebra structure on the direct sum of the underlying vector space and
a finite-dimensional vector space. In addition, their matched pair is constructed and related
to the matched pair of their sub-adjacent Lie algebras. Besides, Zinbiel coalgebras are
introduced, and linked to their underlying Lie coalgebras and coassociative algebras.
Moreover, the related Manin triple is defined and used to characterize Zinbiel bialgebras,
and their equivalence to the associated matched pair.

\section{Zinbiel algebras: basic properties and consequences}
\subsection{Zinbiel algebras}
\begin{definition}
A left (resp. right) Zinbiel algebra is a vector space $\mathcal{A}$ endowed with
a bilinear product $\ast$ satisfying, for all $x,y,z\in \mathcal{A}$,
\begin{eqnarray}\label{eq_left_Zinbiel}
(x\ast y)\ast z=x\ast (y\ast z)+x\ast (z\ast y),
\end{eqnarray}
\begin{eqnarray}\label{eq_right_Zinbiel}
(\quad \mbox{resp.} \quad\; x\ast (y\ast z)=(x\ast y)\ast z+(y\ast x )\ast z\;\quad),
\end{eqnarray}
or, equivalently,
	\begin{eqnarray}
(x,y,z)=x\ast (z\ast y)
	\end{eqnarray}
	\begin{eqnarray}
	(\quad \mbox{resp.} \quad\;  (x,y,z)=-(y\ast x)\ast z)\quad),
	\end{eqnarray}
 where,  $\forall x,y, z\in \A, (x,y,x)=(x\ast y)\ast z-x\ast (y\ast z)$ is the associator associated  to $\ast$.
\end{definition}
\begin{proposition}
	Let
	$\mathcal{A}$ be a vector space,
	$\mu:\mathcal{A}\times \mathcal{A}\rightarrow \mathcal{A}$, 
	$(x,y)\mapsto \mu(x,y)=x\ast y,$ a bilinear product,
	$\tau$   an exchange map defined on 
	$\mathcal{A}\otimes\mathcal{A},$ and $\id$ the identity map defined on $\mathcal{A}$. Then,
		if 
		$(\mathcal{A}, \ast)$ is a left Zinbiel algebra,  the   relation 
		\begin{eqnarray}\label{eq_relation_left_Zinbiel}
	x\ast(y\ast z)=y\ast(x\ast z),\quad \forall x,y, z\in \mathcal{A},
		\end{eqnarray}
		equivalent to
		the two  equations: 
		\begin{eqnarray}\label{eq_relation_left1_Zinbiel}
		\mu\circ(\id\otimes \mu)=\mu\circ(\id\otimes\mu)\circ(\tau\otimes\id),
		\end{eqnarray}
	and 
	\begin{eqnarray}\label{eq_relation_left2_Zinbiel}
		\mu\circ(\id\otimes \mu)=(\mu\circ\tau)\circ(\mu\otimes\id)\circ(\id\otimes\tau),
	\end{eqnarray}
	holds.
	
Similarly, if 
$(\mathcal{A}, \ast)$ is a right Zinbiel algebra,  the  relation
\begin{eqnarray}\label{eq_relation_right_Zinbiel}
(x\ast y)\ast z=(x\ast z)\ast y,\quad  \forall x, y, z\in \mathcal{A},
\end{eqnarray}
equivalent to the two  identities:
\begin{eqnarray}\label{eq_relation_right1_Zinbiel}
\mu\circ(\mu\otimes \id)=\mu\circ(\mu\otimes \id)\circ(\id\otimes\tau),
\end{eqnarray}
and 
\begin{eqnarray}\label{eq_relation_right2_Zinbiel}
\mu\circ(\mu\otimes \id )=(\mu\circ\tau)\circ (\id\otimes\mu)\circ(\tau\otimes\id).
\end{eqnarray}
is satisfied.
\end{proposition}
\subsection{Basic properties}

\begin{definition}
	The {opposite} algebra of the algebra $ (\mathcal{A}, \ast )$  is the algebra denoted by
	$ \mathcal{A}^{opp}=(\mathcal{A}, \ast_{opp})$  whose product is given by 
	$ x\ast_{opp}y=y\ast x$, for all $ x, y\in \mathcal{A}$.
\end{definition}
\begin{remark}
	\begin{itemize}
		\item[(i)] If $ (\mathcal{A}, \ast)$ is a commutative algebra, then $ \mathcal{A}^{opp}=(\mathcal{A}, \ast)$. Conversely, 
		if $ \mathcal{A}^{opp}=(\mathcal{A}, \ast)$, then $ (\mathcal{A}, \ast)$ is a commutative algebra.
		\item[(ii)] The { opposite algebra} of the left (resp. right) Zinbiel algebra is a right (resp.  left) Zinbiel algebra under the 
		same underlying vector space.
	\end{itemize}
\end{remark}
In the sequel,
	both the 
	left and 
	right 
	Zinbiel algebras are simply called 
Zinbiel algebras.
\begin{proposition}\cite{Aguiar}
	Let $(\mathcal{A}, \ast)$ be a Zinbiel algebra. Then $(\mathcal{A}, \{,\}_{\ast}),$
	such that,
	for all $x,y\in \mathcal{A}$, $\{x,y\}_{\ast}=x\ast y+y\ast x,$
	is a commutative associative algebra.
\end{proposition}
\begin{proposition}
Consider a Zinbiel algebra $(\mathcal{A}, \mu),$ $\mu:\mathcal{A}\times \mathcal{A}\rightarrow \mathcal{A}$, 
$(x,y)\mapsto \mu(x,y)=x\ast y,$ and the bilinear  map $\tau:\mathcal{A}\otimes\mathcal{A}\rightarrow\mathcal{A}\otimes\mathcal{A}$ given by 
$\tau(x\otimes y)=y\otimes x$ for all $x, y \in \mathcal{A}$. Then,
the following identities hold:
	\begin{eqnarray}\label{eq_useful1}
	\mu\circ(\id\otimes(\mu\circ\tau))=\mu\circ(\mu\otimes\id)\circ(\id\otimes\tau)+
	(\mu\circ\tau)\circ(\id\otimes(\mu\circ\tau))\circ(\tau\otimes\id)
	\end{eqnarray}
	\begin{eqnarray}\label{eq_useful2}
	(\mu\circ\tau)\circ(\mu\otimes\id)=\mu\circ(\mu\otimes\id)\circ(\id\otimes\tau)+
	(\mu\circ\tau)\circ(\id\otimes(\mu\circ\tau))\circ(\tau\otimes\id)
	\end{eqnarray}
		\begin{eqnarray}\label{eq_useful3}
		(\mu\circ\tau)\circ((\mu\circ\tau)\otimes\id)=(\mu\circ\tau)\circ(\id\otimes\mu)+
	(\mu\circ\tau)\circ(\id\otimes(\mu\circ\tau))
	\end{eqnarray}
	\begin{eqnarray}\label{eq_useful4}
	\mu\circ(\id\otimes\mu)=\mu\circ(\id\otimes\mu)\circ(\tau\otimes\id)
	\end{eqnarray}
\end{proposition}
\begin{theorem}\label{theo_comm_Zinbiel}
	Let $(\A, \ast)$ be a Zinbiel algebra. Then $(\A, [,]_{\ast})$ is a Lie algebra,  where,  for all $x,y\in \A,$
	$[x,y]_{\ast}=x\ast y -y\ast x$. 
\end{theorem}
\textbf{Proof.}

Consider a Zinbiel algebra $(\A, \mu)$,  where for all $x,y\in \A, \,\mu(x,y)=x \ast y$.
By its definition, the commutator of a bilinear product is bilinear and skew symmetric. Thus, proving that 
$(\A, [,]:=\mu-\mu\circ\tau)$ is a Lie algebra, it only  remains to prove that $[,]:=\mu-\mu\circ\tau$
satisfies the Jacobi  identity. Indeed,  we have, for all $x,y,z\in \A$, 
\beqs
[x,[y,z]]+[y,[z,x]]+[z, [x,y]]&=&[x, y\ast z-z\ast y]+[y, z\ast x-x\ast z]+[z, x\ast y-x\ast y]\cr 
&=&\{x\ast (y\ast z)-y\ast(x\ast z)\}+\{y\ast(z\ast x)-z\ast(y\ast x)\}\cr 
&+&\{z\ast(x\ast y)-x\ast(z\ast y)\}+\{(y\ast x)\ast z-(y\ast z)\ast x\}\cr 
&+&\{(z\ast y)\ast x-(z\ast x)\ast y\}+\{(x\ast z)\ast y-(x\ast y)\ast z\}\cr 
&=&
\{\mu\circ(\id\otimes\mu)-\mu\circ(\id\otimes\mu)\circ(\tau\otimes\id)\}(x,y,z)\cr
&+&
\{\mu\circ(\id\otimes\mu)-\mu\circ(\id\otimes\mu)\circ(\tau\otimes\id)\}(y,z,x)\cr
&+&
\{\mu\circ(\id\otimes\mu)-\mu\circ(\id\otimes\mu)\circ(\tau\otimes\id)\}(z,x,y)\cr 
&+&
\{(\mu\circ\tau)\circ(\id\otimes(\mu\circ\tau))\cr &-&\mu\circ((\mu\circ\tau)\otimes\id)\circ(\id\otimes\tau)\}(z,x,y)\cr
&+&
\{(\mu\circ\tau)\circ(\id\otimes(\mu\circ\tau))\cr 
&-&\mu\circ((\mu\circ\tau)\otimes\id)\circ(\id\otimes\tau)\}(x,y,z)\cr 
&+&
\{(\mu\circ\tau)\circ(\id\otimes(\mu\circ\tau))\cr 
&-&\mu\circ((\mu\circ\tau)\otimes\id)\circ(\id\otimes\tau)\}(y,z,x).
\eeqs
Using 
 the relations \eqref{eq_relation_left1_Zinbiel}, \eqref{eq_relation_left2_Zinbiel}, 
\eqref{eq_relation_right2_Zinbiel} and \eqref{eq_relation_right1_Zinbiel}, 
the right hand side of the last  equality of this equation vanishes. Therefore, for all $x, y, z\in \A$, 
\beqs
[x,[y,z]_{\ast}]_{\ast}+[y,[z,x]_{\ast}]_{\ast}+[z, [x,y]_{\ast}]_{\ast}=0,
\eeqs
i.e.  
the  Jacobi identity associated to the bilinear product $[,]_{\ast}$ holds. Hence,  $(\A, [,]_{\ast})$ is a Lie algebra.
$\cqfd$

From the definition of
a center-symmetric algebra given in \cite{Hounkonnou_D_CSA} by,  
 for all
$x,y,z\in \mathcal{A},$
$(x,y,z)_{\circ}:=(z,y,x)_{\circ},$
where   $(x,y,z)_{\circ}:=(x\circ y)\circ z-x\circ(y\circ z)$	is the associator of 
 the bilinear product $\circ$, we have:
	
\begin{proposition}\label{theo_Zinbiel_center_symmetric}
	Any Zinbiel algebra is a center-symmetric algebra.
\end{proposition}
\begin{example}
Consider a vector space $\mathcal{A}$ given by $\displaystyle \mathcal{A}=\mathbb{C}[X]:=\left\lbrace a=\sum_{k\in \mathbb{N}} a_k X^k, a_k\in \mathbb{C}\right\rbrace$.
On $\mathcal{A}$, we define the bilinear products given by, for all 
$\displaystyle a, b\in \mathcal{A},$
\begin{eqnarray}
	a\ast b:= b\int_{0}^{X} a(t) dt\quad \mbox{ and }\quad
\displaystyle a \circ b:=\int_{0}^{X}b(t)\partial_{t}(a(t)) dt,
\end{eqnarray}
where 
$\partial_t:=\frac{d}{d_t}$.
Then, the associators of these products  are given,  for any $a, b, c\in \mathcal{A},$ by 
\begin{eqnarray}
(a,b,c)_{\ast }&=&-(b\ast a)\ast c,\\
(a,b,c)_{\circ}&=& a\circ (c\circ b).
\end{eqnarray} 
Therefore, $(\mathbb{C}[X], \ast)$ and
$(\mathbb{C}[X], \circ)$ are  right and left Zinbiel algebras, respectively.
\end{example}
\subsection{Bimodule}
\begin{definition}
	A bimodule of a Zinbiel algebra is a triple $(l,r, \mathcal{A}),$ where 
	$\mathcal{A}$ is a vector space endowed with a  Zinbiel algebra  structure,  
	$\displaystyle V$ is a vector space, and   
	$ \displaystyle l,r : \mathcal{A} \rightarrow \mathfrak{gl}(V)$ are two linear maps satisfying the following relations, 
	for all $ x, y \in \mathcal{A},$
	\begin{eqnarray}\label{eq_bimod_Zinbiel_1}
		l_xl_y=l_{x\cdot y}+l_{y\cdot x},
	\end{eqnarray}
	\begin{eqnarray}\label{eq_bimod_Zinbiel_2}
		l_xr_y=r_{x\cdot y}=r_{y}r_{x}+r_{y}l_x.
	\end{eqnarray}  
\end{definition}
\begin{proposition} 
	Let $(\mathcal{A}, \cdot)$ be a Zinbiel algebra. Consider a vector space $V$  over a field $\mathcal{K}$
	and two linear maps
	$\displaystyle l, r : \mathcal{A}\rightarrow \mathfrak{gl}(V)$. The triple $(l, r, V)$ is a bimodule of $\mathcal{A}$ if and only if
	there is a Zinbiel algebra structure on the vector space $\mathcal{A}\oplus V$ given by, $\ \forall x, y \in \mathcal{A}$ and all $u, v \in V,$
	\begin{eqnarray}\label{eq_Zinbiel_bimodule}
		(x+u)\ast (y+v) = x\cdot y+(l_{x}v+r_{y}u).
	\end{eqnarray} 
\end{proposition}

\begin{proposition}
	Let $(l,r,V)$ be a bimodule of a Zinbiel algebra $(\mathcal{A}, \cdot)$, where $V$ is a vector space and 
	$l,r:\mathcal{A}\rightarrow \mathfrak{gl}(V)$ are two linear maps. Then,
	\begin{enumerate}
		\item the  following conditions are satisfied, 
		for all $x,y\in \mathcal{A}$,
		\begin{eqnarray}\label{eq_relation_suple_bimodule_Zinbiel_1}
			l_{xy}=r_yl_x,
		\end{eqnarray}
		\begin{eqnarray}\label{eq_relation_suple_bimodule_Zinbiel_2}
			r_xr_y=r_yr_x.
		\end{eqnarray}
		\item  the linear map 
		$l-r: \mathcal{A} \rightarrow \mathfrak{gl}(V), x \mapsto l_x-r_x$
		is a  representation of the  sub-adjacent
		Lie  algebra $\mathcal{G}(\mathcal{A}):=(\mathcal{A}, [,]_{\cdot})$ associated to $(\mathcal{A}, \cdot)$.
	\end{enumerate}
\end{proposition}
\subsection{Matched pair}
\begin{definition}\cite{Ni_B}\label{bimodcommass} 
	Let $\displaystyle (\mathcal{G}, \cdot)$ and $\displaystyle (\mathcal{H}, \circ)$ be two commutative associative  algebras, and    
	${ \displaystyle \rho: \mathcal{H} \rightarrow \mathfrak{gl}(\mathcal{G})}$ and 
	${ \displaystyle \mu: \mathcal{G}\rightarrow \mathfrak{gl}( \mathcal{H})}$ be two ${ \mathcal{K}}$-linear maps which are    
	representations of ${ \mathcal{H}}$ and ${ \mathcal{G}}$, respectively,  satisfying the following relations:
	$ \mbox{for all } x, y \in \mathcal{G}$ and all ${ \displaystyle a, b \in~\mathcal{H}},$ 
	\begin{eqnarray}\label{eq_matched_asso_commu_1}
	{ \mu(x)(a\circ b)=(\mu(x)a)\circ b+\mu(\rho(a)x)b},
	\end{eqnarray}
	\begin{eqnarray}\label{eq_matched_asso_commu_2}
	{ \rho(a)(x\cdot y)=(\rho(a)x)\cdot b+\rho(\mu(x)a)y}.
	\end{eqnarray}
	Then, ${ \displaystyle(\mathcal{G}, \mathcal{H}, \rho, \mu)}$ is called a matched pair of the commutative associative algebras
	${ \displaystyle \mathcal{G}}$ and ${ \displaystyle \mathcal{H}},$ 
	denoted by ${ \displaystyle \mathcal{G}\bowtie^{\mu}_{\rho}\mathcal{H}}.$
	
	In this case,  ${ \displaystyle (\mathcal{G}\oplus \mathcal{H}, \ast)}$ defines a commutative  associative algebra with respect 
	to the product ${ \ast},$  given, for all ${ x,y\in \mathcal{G}}$ and all ${ a,b\in \mathcal{H}}$, by
	\begin{center}
		${ \displaystyle (x+a)\ast (y+b)= x\cdot y+\mu(a)y+\mu(b)x+a\circ b+\rho(x)b+\rho(y)a}$.
	\end{center} 
\end{definition}
\begin{definition}\cite{Madjid} 
	Let ${ \displaystyle (\mathcal{G}, [,]_{_\mathcal{G}})}$ and 
	${ \displaystyle (\mathcal{H}, [,]_{_\mathcal{H}})}$ be two Lie algebras and let  
	${ \displaystyle \mu: \mathcal{H} \rightarrow \mathfrak{gl}(\mathcal{G})}$ and 
	${ \displaystyle \rho: \mathcal{G}\rightarrow \mathfrak{gl}( \mathcal{H})}$ be two Lie algebra representations satisfying the following relations, for all 
	${ \displaystyle x, y \in \mathcal{G}}$ and all $a, b \in~\mathcal{H},$ 
	\begin{eqnarray}\label{eq_matched_Lie_commu_1}
	{ \rho(x)\left[a, b\right]_{_\mathcal{H}}-\left[\rho(x)a, b\right]_{_\mathcal{H}}
		-\left[a, \rho(x)b\right]_{_\mathcal{H}}+\rho(\mu(a)x)b-\rho(\mu(b)x)a=0},
	\end{eqnarray}
	\begin{eqnarray}\label{eq_matched_Lie_commu_2}
	{ \mu(a)\left[x, y\right]_{_\mathcal{G}}-\left[\mu(a)x, y\right]_{_\mathcal{G}}
		-\left[x, \mu(a)y\right]_{_\mathcal{G}}+\mu(\rho(x)a)y-\mu(\rho(y)a)x=0}.
	\end{eqnarray}
	Then, ${ \displaystyle(\mathcal{G}, \mathcal{H}, \rho, \mu)}$ is called a matched pair of the Lie algebras
	${ \displaystyle \mathcal{G}}$ and ${ \displaystyle \mathcal{H}},$ 
	denoted by ${ \displaystyle \mathcal{H}\bowtie_{\mu}^{\rho}\mathcal{G}}.$
	In this case,  ${ \displaystyle (\mathcal{G}\oplus \mathcal{H}, [,]_{_{\mathcal{G}\oplus\mathcal{H}}})}$ is a Lie algebra with respect to the product
	${ [,]_{_{\mathcal{G}\oplus\mathcal{H}}}}$ defined on the direct sum vector space ${ \mathcal{G}\oplus\mathcal{H}}$  by, 
	${ \forall x,y\in \mathcal{G}}$ and all 
	${ a,b\in \mathcal{H}}$, 
	\begin{eqnarray}
	{ \displaystyle [(x+a), (y+b)]_{_{\mathcal{G}\oplus\mathcal{H}}}= [x, y]_{_\mathcal{G}}+\mu(a)y-\mu(b)x+[a, b]_{_\mathcal{H}}+\rho(x)b-\rho(y)a}.
	\end{eqnarray} 
\end{definition}
\begin{theorem} 
	Let  ${ (\mathcal{A}, \cdot)}$ and ${ (\mathcal{B}, \circ)}$  be two Zinbiel algebras. Consider the four linear maps defined as 
	${ l_{\mathcal{A}}, r_{\mathcal{A}}: \mathcal{A}\rightarrow \mathfrak{gl}(\mathcal{B})}$ 
	and ${ l_{\mathcal{B}}, r_{\mathcal{B}}: \mathcal{B}\rightarrow \mathfrak{gl}(\mathcal{A})}$. 
	Suppose that ${ (l_{\mathcal{A}}, r_{\mathcal{A}}, \mathcal{B})}$ and ${ (l_{\mathcal{B}}, r_{\mathcal{B}},\mathcal{A})}$
	are bimodules of ${ \mathcal{A}}$ and ${ \mathcal{B}}$, respectively,   obeying the relations,
	for all ${ x, y \in \mathcal{A}}$ and all ${ a, b \in \mathcal{B}}$, 
	\begin{eqnarray}\label{eqmatchzin1}
	{  r_{\mathcal{B}}(a)(x\cdot y+ y \cdot x)-x\cdot (r_{\mathcal{B}}(a)y)-r_{\mathcal{B}}(l_{\mathcal{A}}(y)a)x=0},  
	\end{eqnarray}
	\begin{eqnarray}\label{eqmatchzin2}
	{ r_{\mathcal{A}}(x)(a\circ b+ b \circ a)-a\circ (r_{\mathcal{A}}(x)b)-r_{\mathcal{A}}(l_{\mathcal{B}}(b)x)a=0}, 
	\end{eqnarray}
	\begin{eqnarray}\label{eqmatchzin3}
	{  
		l_{\mathcal{B}}(a)(x\cdot y)}
	&=&
	{  ((l_{\mathcal{B}}+r_{\mathcal{B}})(a)x)\cdot y+l_{\mathcal{B}}((l_{\mathcal{A}}+r_{\mathcal{A}})(x)a)y}\\
	&=&
	{ 
		x\cdot (l_{\mathcal{B}}(a)y)+r_{\mathcal{B}}(r_{A}(y)a)x},\nonumber
	\end{eqnarray}
\begin{eqnarray}\label{eqmatchzin4}
{  l_{\mathcal{A}}(x)(a\circ b)}&=&
{ l_{\mathcal{A}}((l_{\mathcal{B}}+r_{\mathcal{B}})(a)x)b+
	((l_{\mathcal{A}}+r_{\mathcal{A}})(x)a) \circ b}\\
&=&
{  a\circ (l_{\mathcal{A}}(x)b)+r_{\mathcal{A}}(r_{\mathcal{B}}(b)x)a}.\nonumber
\end{eqnarray}
Then, there is a Zinbiel algebra structure on ${ \mathcal{A}\oplus \mathcal{B}}$ 
given by: 
\begin{eqnarray}\label{eqzinbieldouble}
{ \displaystyle (x+a)\ast (y+b)= (x \cdot y + l_{\mathcal{B}}(a)y+r_{\mathcal{B}}(b)x)+ (a \circ b + l_{\mathcal{A}}(x)b+r_{\mathcal{A}}(y)a)}.
\end{eqnarray} 
\end{theorem}

\textbf{Proof.}

Suppose $x, y, z \in \A$ and $a, b, c \in \B.$ Then,  using the bilinear product $\ast$ defined on 
$\A\oplus \B$ by $(x+a)\ast (y+b)=\{x\cdot y+l_{\B}(x)y+r_{\B}(b)x\}+\{a\circ b+l_{\A}(x)b+r_{\A}(y)a\},$ 
 we have:
\beq\label{pap1}
&&\{(x+a)\ast(y+b)\}\ast (z+c) = (x\cdot y)\cdot z+((l_{\B}(a)y)\cdot z +l_{\B}(r_{\A}(y)a z))+l_{\B}(a\circ b)z  \cr 
&&+ r_{\B}(c)(x\cdot y) 
+((r_{\B}(b)x)\cdot z+l_{\B}(l_{\A}(x)b)z)+
r_{\B}(c)(l_{\B}(a)y)+r_{\B}(c)(r_{\B}(b)x) \cr 
&&+ (a\circ b)\circ c  +((l_{\A}(x)b)\circ c+ 
l_{\A}(r_{\B}(b)x)c)+l_{\A}(x\cdot y)c+
r_{\A}(z)(a\circ b)\\
&&+  r_{\A}(z)(l_{\A}(x)b) + ((r_{\A}(y)a)\circ c+  l_{\A}(l_{\B}(a)y)c)+r_{\A}(z)(r_{\A}(y)a), \nonumber
\eeq 
\beq\label{pap2}
&&\{(y+b)\ast(x+a)\}\ast (z+c) = (y\cdot x)\cdot z+((l_{\B}(b)x)\cdot z +l_{\B}(r_{\A}(x)b z))+l_{\B}(b\circ a)z  \cr 
&&+ r_{\B}(c)(y\cdot x) 
+((r_{\B}(a)y)\cdot z+l_{\B}(l_{\A}(y)a)z)+
r_{\B}(c)(l_{\B}(b)x)+r_{\B}(c)(r_{\B}(a)y)  \cr
&&+ (b\circ a)\circ c+((l_{\A}(y)a)\circ c+ 
l_{\A}(r_{\B}(a)y)c)+l_{\A}(y\cdot x)c+
r_{\A}(z)(b\circ a)  \\ 
&&+((r_{\A}(x)b)\circ c+ l_{\A}(l_{\B}(b)x)c) +r_{\A}(z)(l_{\A}(y)a)+r_{\A}(z)(r_{\A}(x)b), \nonumber
\eeq
\beq\label{pap3}
&&(x+a)\ast \{(y+b)\ast (z+c)\}=x\cdot(y\cdot z)+(x\cdot(l_{\B}(b)z)+r_{\B}(r_{\A}(z)b)x)+
l_{\B}(a)(y\cdot z)\cr 
&&+ l_{\B}(a)(l_{\B}(b)z)+l_{\B}(a)(r_{\B}(c)y)+
r_{\B}(b\circ c)x+ (x\cdot(r_{\B}(c)y)+r_{\B}(l_{\A}(y)c)x)\cr
&&+ a\circ (b \circ c) + (a\circ (l_{\A}(y)c)+r_{\A}(r_{\B}(c)y)a)+
(a\circ (r_{\A}(z)b)+r_{\A}(l_{\B}(b)z)a)
\\&&+ l_{\A}(x)(b\circ c)
 +l_{\A}(x)(l_{\A}(y)c)+l_{\A}(x)(r_{\A}(z)b)+r_{\A}(y\cdot z)a. \nonumber
\eeq
Summing  \eqref{pap1} and \eqref{pap2} gives \eqref{pap3}, what is  equivalent to say that
$
(l_{\A}, r_{\A}, \B)$ and  $(l_{\B}, r_{\B}, \A)$ are bimodules of the Zinbiel algebras 
$(\A, \cdot )$  and $( \B, \circ)$,  respectively, and the four linear maps 
$l_{\A}, r_{\A}, l_{\B} $,  and  $r_{\B}$ satisfy   the equations  \eqref{eqmatchzin1}, 
\eqref{eqmatchzin2}, \eqref{eqmatchzin3}, and  \eqref{eqmatchzin4}.
$\cqfd$

\begin{corollary}
Let ${ (\mathcal{A}, \mathcal{B}, l_{\mathcal{A}}, r_{\mathcal{A}}, l_{\mathcal{B}}, r_{\mathcal{B}})}$ be a matched pair of the Zinbiel algebras 
${ (\mathcal{A}, \cdot)}$ and ${ (\mathcal{B}, \circ )}.$ Then,
\begin{enumerate}
	\item
	${ (\mathcal{G}_{\mbox{ass}}(\mathcal{A}), \mathcal{G}_{\mbox{ass}}(\mathcal{B}), l_{\mathcal{A}}+r_{\mathcal{A}}, l_{\mathcal{B}}+r_{\mathcal{B}})}$ 
	is a matched pair of the   commutative associative algebras 
	${ \mathcal{G}_{\mbox{ass}}(\mathcal{A}):=(\mathcal{A}, \{,\}_{_\cdot})}$ and  
	${ \mathcal{G}_{\mbox{ass}}(\mathcal{B}):=(\mathcal{B}, \{,\}_{_\ast})}.$
	\item
	${ (\mathcal{G} (\mathcal{A}), \mathcal{G} (\mathcal{B}), l_{\mathcal{A}}-r_{\mathcal{A}}, l_{\mathcal{B}}-r_{\mathcal{B}})}$ 
	is a matched pair of the   Lie algebras ${ \mathcal{G}(\mathcal{A}):=(\mathcal{A}, [,]_{_\cdot})}$ and  ${ \mathcal{G}(\mathcal{B}):=(\mathcal{B},[,]_{_\circ})}.$
\end{enumerate}
\end{corollary}
\section{Zinbiel coalgebras, Manin triple and bialgebra}
 The structures of nonassociative coalgebras  and their basic properties are investigated and discussed in \cite{Anquela_C} 
	and references therein. This section is devoted to the construction of Zinbiel bialgebra using their associated coalgebras and 
	Manin triple.
\subsection{Zinbiel coalgebras}
\begin{definition}
	Let ${ \mathcal{A}}$ be a vector space equipped with the linear map 
	${ \Delta :\mathcal{A}\rightarrow \mathcal{A}\otimes \mathcal{A}}$. 
	The couple ${ (\mathcal{A}, \Delta)}$ is a { right Zinbiel coalgebra}  if 
	${ \Delta}$ satisfies the following identity
	\begin{eqnarray}\label{eq_right_Zinbiel_coalgebra}
	{ 
		(\id\otimes\Delta)\circ\Delta =
		(\Delta\otimes\id)\circ\Delta+((\tau\circ\Delta)\otimes\id)\circ\Delta
	},
	\end{eqnarray}
	illustrated by the following commutative diagram:
	\begin{eqnarray}
	\SelectTips{cm}{10}
	\xymatrix{
		{ \mathcal{A}}\ar[rr]^-{{ \Delta}} \ar[d]_-{{ \Delta}}& &
		{ \mathcal{A} \otimes \mathcal{A}} \ar[d]^-{{ \id\otimes \Delta }}\\
		{ \mathcal{A}\otimes\mathcal{A}} \ar[rr]_-{{ (\Delta+(\tau\circ\Delta))\otimes\id}}
		& &  { \mathcal{A} \otimes \mathcal{A}\otimes \mathcal{A}}
	}
	\end{eqnarray}	
	where ${ \id}$ is the identity map on ${ \mathcal{A}},$ and ${ \tau}$ is the exchange map defined on ${ \mathcal{A} \otimes \mathcal{A}}$.
\end{definition}
\begin{definition} 
Let ${ \mathcal{A}}$ be a vector space endowed  with the linear map ${ \Delta :\mathcal{A}\rightarrow \mathcal{A}\otimes \mathcal{A}}$. 
The couple ${ (\mathcal{A}, \Delta)}$ is a { left Zinbiel coalgebra}  if ${ \Delta}$ satisfies the following identity: 
\begin{eqnarray}\label{eq_left_Zinbiel_coalgebra}
{ 
	(\Delta\otimes\id)\circ\Delta =
	(\id\otimes\Delta)\circ\Delta+(\id\otimes (\tau\circ\Delta))\circ\Delta},
\end{eqnarray}
illustrated by the following commutative diagram:
\begin{eqnarray}
\SelectTips{cm}{10}
\xymatrix{
	{ \mathcal{A}} \ar[rr]^-{{ \Delta}} \ar[d]_-{{ \Delta}}& &
	{ \mathcal{A} \otimes \mathcal{A}} \ar[d]^-{{ \Delta\otimes \id}}\\
	{ \mathcal{A}\otimes\mathcal{A}} \ar[rr]_-{{ \id\otimes (\Delta+(\tau\circ\Delta))}} &&
	{ \mathcal{A} \otimes \mathcal{A} \otimes\mathcal{A}}
}
\end{eqnarray}
where ${ \id}$ is the identity map on ${ \mathcal{A}},$ and ${ \tau}$ is the exchange map defined on ${ \mathcal{A} \otimes \mathcal{A}}$.
\end{definition}
\begin{proposition}
Let ${ \mathcal{A}}$ be a vector space equipped with a linear map 
${ \Delta: \mathcal{A}\rightarrow \mathcal{A}\otimes\mathcal{A}}$,
${ \tau}$ be the exchange map defined on ${ \mathcal{A}\otimes\mathcal{A}},$ 
and ${ \id}$ be the identity map defined on ${ \mathcal{A}}$.
\begin{enumerate}
\item If ${ (\mathcal{A}, \Delta)}$ is a right Zinbiel coalgebra, then the relation 
\begin{eqnarray}\label{eq_relation_right_coalgebra}
{ 
	(\id\otimes\Delta)\circ\Delta=
	(\tau\otimes\id)\circ(\id\otimes\Delta)\circ\Delta
	=(\tau\otimes\id)\circ(\Delta\otimes\id)\circ(\tau\circ\Delta)}
\end{eqnarray}	
is satisfied.
\item 
If ${ (\mathcal{A}, \Delta)}$ is a left Zinbiel coalgebra, then the  relation 
\begin{eqnarray}\label{eq_relation_left_coalgebra}
{ 
	(\Delta\otimes\id)\circ\Delta=
	(\id\otimes\tau)\circ(\Delta\otimes\id)\circ\Delta=
	(\id\otimes\tau)\circ(\id\otimes\Delta)\circ(\tau\circ\Delta)}
\end{eqnarray}
holds.
\end{enumerate} 
\end{proposition}
\begin{definition} 
A commutative infinitesimal bialgebra ${ (\mathcal{A}, \cdot, \Delta)}$ is an infinitesimal bialgebra 
${ (\mathcal{A}, \cdot, \Delta)}$ such that ${ (\mathcal{A}, \cdot)}$ is a commutative algebra and 
${ (\mathcal{A}, \Delta)}$ 
is a cocommutative coalgebra.
\end{definition}
\begin{definition} 
The opposite coalgebra of the coalgebra ${ (\mathcal{A}, \Delta)}$ is a coalgebra 
${ (\mathcal{A}, \tau\circ\Delta)},$ where 	the coproduct ${ \tau\circ\Delta}$ is defined by:	
${ \tau\circ\Delta: \mathcal{A}\longrightarrow \mathcal{A}\otimes\mathcal{A}}$.
\end{definition}
\begin{proposition} 
The { opposite} coalgebra of a { right Zinbiel} coalgebra is a
{ left Zinbiel} coalgebra. Conversely, 
the { opposite} coalgebra of a { left Zinbiel} coalgebra is a 
{right Zinbiel} coalgebra.
\end{proposition}
\begin{proposition}
Consider a right Zinbiel coalgebra ${ (\mathcal{A}, \Delta)}$, and ${ \tau}$  the exchange map  given by 
${ \tau:\mathcal{A}\otimes\mathcal{A}\rightarrow\mathcal{A}\otimes\mathcal{A}}$,
for all ${ x,y\in \mathcal{A}}$,
${ \tau(x\otimes y)=y\otimes x}$.
Then, the following equations hold:
\begin{eqnarray}\label{eq_useful_dual1}
{ 
(\id\otimes(\tau\circ\Delta))\circ\Delta=(\id\otimes\tau)\circ(\Delta\otimes\id)\circ\Delta+
(\tau\otimes\id)\circ(\id\otimes(\tau\circ\Delta))\circ(\tau\circ\Delta)},
\end{eqnarray}
\begin{eqnarray}\label{eq_useful_dual2}
{ 
(\Delta\otimes\id)\circ(\tau\circ\Delta)=(\id\otimes\tau)\circ(\Delta\otimes\id)\circ\Delta+
(\tau\otimes\id)\circ(\id\otimes(\tau\circ\Delta)\circ(\tau\circ\Delta)},
\end{eqnarray}
\begin{eqnarray}\label{eq_useful_dual3}
{ 
((\tau\circ\Delta)\otimes\id)\circ(\tau\circ\Delta)=
(\id\otimes\Delta)\circ(\tau\circ\Delta)+
(\id\otimes(\tau\circ\Delta))\circ(\tau\circ\Delta)}.
\end{eqnarray} 
\end{proposition}

\textbf{Proof.}

Consider a right Zinbiel coalgebra $(\A, \Delta)$ and  the exchange map $\tau$  on $\A\otimes\A$. We have:
\beqs
(\id\otimes(\tau\circ\Delta))\circ\Delta&=&
(\id\otimes\tau)\circ(\id\otimes\Delta)\circ\Delta\cr 
&=&(\id\otimes\tau)\circ((\Delta\otimes\id)\circ\Delta)+
(\id\otimes\Delta)\circ((\tau\otimes\Delta)\otimes\id)\circ\Delta\cr 
&=&(\id\otimes\Delta)\circ(\Delta\otimes\id)\circ\Delta+
(\id\otimes\tau)\circ(\id\otimes(\tau\circ\Delta))\circ(\tau\circ\Delta)\cr 
&=&(\id\otimes\tau)\circ(\Delta\otimes\id)\circ\Delta+(\id\otimes\Delta)\circ(\tau\circ\Delta)\cr 
&=&(\id\otimes\tau)\circ(\Delta\otimes\id)\circ\Delta+(\Delta\otimes\id)\circ\Delta\cr 
&=&(\id\otimes\tau)\circ(\Delta\otimes\id)\circ\Delta+
(\tau\otimes\id)\circ((\tau\circ\Delta)\otimes\id)\circ\Delta\cr 
(\id\otimes(\tau\circ\Delta))\circ\Delta
&=& (\id\otimes\tau)\circ(\Delta\otimes\id)\circ\Delta+
(\tau\otimes\id)\circ(\id\otimes(\tau\circ\Delta))\circ(\tau\circ\Delta).
\eeqs
Hence,  the relation \eqref{eq_useful_dual1} holds.

Besides, 
\beqs
(\id\otimes(\tau\circ\Delta))\circ\Delta&=&
((\tau\circ\Delta)\otimes\id)\circ(\tau\circ\Delta)
=(\tau\otimes)\circ(\Delta\otimes\id)\circ(\tau\circ\Delta)\cr 
&=&(\id\otimes\Delta)\circ\Delta=
(\Delta\otimes\id)\circ(\tau\circ\Delta),
\eeqs
which is the relation \eqref{eq_useful_dual2}.

Using  the following identities,
\beqs
((\tau\circ\Delta)\otimes\id)\circ(\tau\circ\Delta)&=&
(\id\otimes(\tau\circ\Delta))\circ\Delta,\cr
(\id\otimes\Delta)\circ(\tau\circ\Delta)&=&
(\Delta\otimes\id)\circ\Delta=
(\tau\otimes\id)\circ((\tau\circ\Delta)\otimes\id)\circ\Delta\cr 
(\id\otimes\Delta)\circ(\tau\circ\Delta)&=&
(\tau\otimes\id)\circ(\id\otimes(\tau\circ\Delta))\circ(\tau\circ\Delta),
\eeqs
and 
\beqs
(\id\otimes(\tau\circ\Delta))\circ(\tau\circ\Delta=
(\id\otimes\tau)\circ (\id\otimes\Delta)\circ(\tau\circ\Delta)=
(\id\otimes\tau)\circ(\Delta\otimes\id)\circ\Delta,
\eeqs
we obtain that the relations \eqref{eq_useful_dual1} and \eqref{eq_useful_dual3} are  equivalent.
	$\cqfd$

\begin{definition}
A cocommutative coassociative coalgebra is the couple  ${ (\mathcal{A},\Delta)}$, where ${ \mathcal{A}}$ is a 
vector space and ${ \Delta: \mathcal{A}\rightarrow \mathcal{A}\otimes\mathcal{A}}$ a linear map such that the following relations hold:
\begin{eqnarray}\label{eq_cocommutator_Zinbiel}
{ 
\Delta=\tau\circ\Delta},
\end{eqnarray} 
\begin{eqnarray}\label{eq_coassociator_Zinbiel}
{ 
(\Delta\otimes\id)\circ\Delta=(\id\otimes\Delta)\circ\Delta}.
\end{eqnarray} 
\end{definition}
\begin{proposition}
Consider  a vector space  ${ \mathcal{A}}$  and ${ \Delta}$ a linear map given by ${ \Delta: \mathcal{A} \rightarrow \mathcal{A}\otimes\mathcal{A}}$. 
The couple ${ (\mathcal{A}, \Delta)}$ is a Zinbiel coalgebra if  ${ (\mathcal{A}, \Delta_{\{,\}})}$ is a cocommutative coassociative 
coalgebra, where ${ \Delta_{\{,\}}: \mathcal{A} \rightarrow \mathcal{A}\otimes \mathcal{A}}$ is given by ${ \Delta_{\{,\}}:=\Delta+\tau\circ\Delta},$ 
and ${ \tau}$ is the exchange map defined on ${ \mathcal{A}\otimes\mathcal{A}}$.
\end{proposition}
\begin{definition}
Let ${ \mathcal{A}}$ be a vector space equipped with the linear map
${ \Delta: \mathcal{A}\rightarrow \mathcal{A}\otimes\mathcal{A}}$. 
Then, the couple ${ (\mathcal{A}, \Delta)}$ is a Lie coalgebra if 
the following relations are satisfied:
\begin{eqnarray}\label{eq_skew_Lie_coalgebra}
{ 
\Delta=-\tau\circ\Delta},
\end{eqnarray}
\begin{eqnarray}\label{eq_Jacobi_Lie_coalgebra}
{ 
(\id\otimes\Delta)\circ\Delta+
(\id\otimes\tau)\circ(\Delta\otimes\id)\circ\Delta-
(\Delta\otimes\id)\circ\Delta=0}.
\end{eqnarray}
\end{definition}
\begin{proposition}
Consider a  vector space  ${ \mathcal{A}}$  equipped with a linear map ${ \Delta: \mathcal{A} \rightarrow \mathcal{A}\otimes\mathcal{A}}$. 
Then, the couple ${ (\mathcal{A}, \Delta)}$ is a Zinbiel coalgebra if the linear map ${ \Delta_{[,]}: \mathcal{A} \rightarrow \mathcal{A}\otimes \mathcal{A}}$ 
given by ${ \Delta_{[,]}:=\Delta-\tau\circ\Delta}$ 
satisfies the  relations \eqref{eq_skew_Lie_coalgebra} and \eqref{eq_Jacobi_Lie_coalgebra}.
\end{proposition}

\textbf{Proof.}

Let $(\A, \mu)$ be a Zinbiel algebra and $(\A, \Delta)$ be a Zinbiel coalgebra. From the definition of 
the exchange map, we have:
\beqs
\Delta_{[,]}&=&\id\circ \Delta_{[,]}=\tau\circ(\tau\circ\Delta_{[,]})=
\tau\circ(\tau\circ(\Delta-\tau\circ\Delta))\cr 
&=&
\tau\circ(\tau\circ\Delta-\Delta)=-\tau\circ(\Delta-\tau\circ\Delta),
\eeqs
and, hence, $\Delta_{[,]}=-\tau\circ\Delta_{[,]}$ which is the relation \eqref{eq_skew_Lie_coalgebra}.

Besides,
\beqs
(\id \otimes \Delta_{[,]} )\circ \Delta_{[,]}&=&
(\id\otimes(\Delta-\tau\circ\Delta))\circ(\Delta-\tau\circ\Delta)\cr
&=&(\id\otimes\Delta)\circ\Delta-(\id\otimes\Delta)\circ(\tau\circ\Delta)-
(\id\otimes(\tau\circ\Delta))\circ\Delta\cr 
&+&(\id\otimes(\tau\circ\Delta))\circ(\tau\circ\Delta).
\eeqs
Then, the   relation 
\beq\label{eq_relation_proof1}
(\id \otimes \Delta_{[,]} )\circ \Delta_{[,]}&=&
(\id\otimes\Delta)\circ\Delta-(\id\otimes\Delta)\circ(\tau\circ\Delta)-
(\id\otimes(\tau\circ\Delta))\circ\Delta\cr 
&+&(\id\otimes(\tau\circ\Delta))\circ(\tau\circ\Delta)
\eeq
holds. In addition, 
\beqs
(\Delta_{[,]}\otimes\id)\circ\Delta_{[,]}
&=&((\Delta-\tau\circ\Delta)\otimes\id)\circ(\Delta-\tau\circ\Delta)\cr 
&=&(\Delta\otimes\id)\circ(\Delta-\tau\circ\Delta))-((\tau\circ\Delta)\otimes\id)\circ(\Delta-(\tau\circ\Delta))
\cr 
&=&(\Delta\otimes\id)\circ\Delta-(\Delta\otimes\id)\circ(\tau\circ\Delta)-
((\tau\circ\Delta)\otimes\id)\circ\Delta\cr
&+&
((\tau\circ\Delta)\otimes\id)\circ(\tau\circ\Delta, 
\eeqs
and the   relation 
\beq\label{eq_relation_proof2}
(\Delta_{[,]}\otimes\id)\circ\Delta_{[,]}
&=&
(\Delta\otimes\id)\circ\Delta-(\Delta\otimes\id)\circ(\tau\circ\Delta)-
((\tau\circ\Delta)\otimes\id)\circ\Delta\cr
&+&
((\tau\circ\Delta)\otimes\id)\circ(\tau\circ\Delta)
\eeq
is satisfied. 

Furthermore, we have 
\beq\label{eq_relation_proof3}
(\id\otimes\tau)\circ(\Delta_{[,]}\otimes\id)\circ\Delta_{[,]}&=&
(\id\otimes\tau)\circ((\tau\circ\Delta)\otimes\id)\circ(\tau\circ\Delta)\cr 
&-&(\id\otimes\tau)\circ(\Delta\otimes\id)\circ(\tau\circ\Delta)
\\&-&
(\id\otimes\tau)\circ((\tau\circ\Delta)\otimes\id)\circ\Delta\cr 
&+&
(\id\otimes\tau)\circ(\Delta\otimes\id)\circ\Delta. \nonumber
\eeq

In addition, by considering the right hand side of the relations \eqref{eq_relation_proof1}, \eqref{eq_relation_proof2},
and \eqref{eq_relation_proof3},  and using the relations  
\eqref{eq_useful_dual1}, \eqref{eq_useful_dual2}, and \eqref{eq_useful_dual3}, we obtain the identity \eqref{eq_Jacobi_Lie_coalgebra}. 
$\cqfd$

\subsection{Manin triple and Zinbiel bialgebras}
\begin{definition}
	Consider a nonassociative algebra ${ (\mathcal{A}, \cdot)}$ and ${ \mathbb{B}}$ a bilinear form on ${ \mathcal{A}}$.
	\begin{itemize}
		\item  ${ \mathbb{B}}$ is symmetric  if ${ \mathbb{B}(x, y)=\mathbb{B}(y, x), \forall x, y \in \mathcal{A}}$.
		\item ${ \mathbb{B}}$ is invariant  if ${ \mathbb{B}(x\cdot y, z)=\mathbb{B}(x, y\cdot z), \forall x, y, z \in \mathcal{A}}$.
		\item  ${ \mathbb{B}}$ is nondegenerate   if the set  ${ \displaystyle \{x \in \mathcal{A}, \mathbb{B}(x, y)=0, \forall y \in \mathcal{A}    \}}$  
		is reduce to ${ \{0\}}$.
	\end{itemize}
\end{definition}
\begin{definition} 
A Manin triple of  Zinbiel algebras is a triple ${ (\mathcal{A}, \mathcal{A}^{+}, \mathcal{A}^{-})}$ 
endowed with a   nondegenerate symmetric bilinear form   ${ \mathbb{B}~(~;~)},$ invariant on ${ \mathcal{A}},$   and satisfying:
\begin{enumerate}
\item ${ \mathcal{A}=\mathcal{A}^{+}\oplus \mathcal{A}^{-}}$ as 
${ \mathbb{K}}$-vector space;
\item ${ \mathcal{A}^{+}}$ and ${ \mathcal{A}^{-}}$ are Zinbiel subalgebras of ${ \mathcal{A}}$; 
\item  ${ \mathcal{A}^{+}}$ and ${ \mathcal{A}^{-}}$ are isotropic with respect  to 
${ \mathbb{B}(;)},$ i.e., 
${ \mathbb{B}(\mathcal{A}^{+};\mathcal{A}^{+})=0= \mathbb{B}(\mathcal{A}^{-};\mathcal{A}^{-})}.$
\end{enumerate}
\end{definition}
\begin{theorem}
Let ${ (\mathcal{A}, \cdot )}$ and ${ (\mathcal{A}^*, \circ )}$ be two Zinbiel algebras. Then, the sixtuple
${ (\mathcal{A}, \mathcal{A}^*, R_{\cdot}^*, L_{\cdot}^*; R_{\circ}^*, L_{\circ}^* )}$ is a matched pair of Zinbiel algebras ${ \mathcal{A}}$ and ${ \mathcal{A}^{*}}$ if and only if   
${ (\mathcal{A}\oplus \mathcal{A}^*, \mathcal{A}, \mathcal{A}^* )}$ is a Manin triple endowed with the bilinear form
${  \mathbb{B}(x+a^*, y+b^*)=<x, b^*>+<y, a^*>,  \forall x, y \in \mathcal{A}}$ and $\forall  a^*, b^* \in \mathcal{A}^*,$ where
${  <,>}$ is the natural pairing between ${  \mathcal{A}}$ and ${ \mathcal{A}^*}.$
\end{theorem}
\begin{proposition}
Let ${ (\mathcal{A}, \cdot)}$ be a Zinbiel algebra and ${ (\mathcal{A}^*, \circ)}$ be a Zinbiel algebra  on its dual space ${ \mathcal{A}^*}.$ Then the following conditions are equivalent:
\begin{enumerate}
\item  ${ (\mathcal{A}\oplus\mathcal{A}^*, \mathcal{A}, \mathcal{A}^*)}$ is the standard Manin triple of considered Zinbiel algebra; 
\item ${ (\mathcal{G}(\mathcal{A}), \mathcal{G}(\mathcal{A}^*),-ad_{\cdot}^*, -ad_{\circ}^*)}$ is a  matched pair of sub-adjacent Lie algebras;
\item  ${ (\mathcal{A}, \mathcal{A}^*, R_{\cdot}^*, L_{\cdot}^*, R_{\circ}^*, L_{\circ}^*)}$ is a matched pair of Zinbiel algebras;
\item ${ (\mathcal{A}, \mathcal{A}^*)}$ is a Zinbiel  bialgebra.
\end{enumerate}
\end{proposition}
\ack	
This work is supported by TWAS Research Grant  RGA No.17-542 RG/ MATHS/AF/AC\_G-FR3240300147. 
The ICMPA-UNESCO Chair is in partnership with the Association pour la 
Promotion Scientifique de l'Afrique (APSA), France, and   Daniel Iagolnitzer Foundation (DIF), 
France, supporting the development of mathematical physics in Africa.


\section*{References}
\begin{thebibliography}{9999}
\bibitem{Santilli}
Santilli R M 1968 Nuovo Cim. Suppl. \textbf{6}  1225-1249
\bibitem{Bai_LSA}
Bai C 2008
Commun. Contemp. Math.,
\textbf{10} No.2 221-260
\bibitem{Loday}
Loday J-L 1995 
Math. Scand. \textbf{77}  189-196
\bibitem{Kawski}
Kawski M 2001
J.  Math.  Sci. \textbf{103} No.6 725-744
\bibitem{Dzhumadil'daev}
Dzhumadil'daev A S 2007 
J.  Math.  Sci., \textbf{144} No.2  3909-3925
\bibitem{Aguiar}
Aguiar M 2000 
Lett. Math. Phys. \textbf{54} 263-277
\bibitem{Hounkonnou_D_CSA}
Hounkonnou M.N.  and Dassoundo  M.L. 2015 Trends Math. 261-273
\bibitem{Ni_B}
Ni X and  Bai C 2013 J. Math. Phys. \textbf{54} 023515
\bibitem{Madjid}
Majid S 1990 Pacific J. Math.
\textbf{141} No.2  311-332
\bibitem{Anquela_C}
Anquela J A,  Cort\'es  T and Montaner F 1994  
\textbf{22} No12  4693-4716
\end{thebibliography}
\end{document}